\documentclass[reqno,12pt]{amsart}

\usepackage[english]{babel} 						
\usepackage[T1]{fontenc}							
\usepackage[utf8]{inputenx}

\usepackage{mathtools}								
\usepackage{empheq}									
\usepackage{enumitem}								

\usepackage{amssymb}								
\usepackage{mathrsfs}								

\usepackage{aliascnt}								
\usepackage{comment}

\usepackage{a4wide,color,eucal}
\usepackage[normalem]{ulem}
\usepackage{epsfig}
\usepackage{a4wide,color,eucal,mathrsfs}
\usepackage[normalem]{ulem}
\usepackage{amsmath,amssymb,epsfig,amsthm} 

\usepackage[babel]{csquotes}							

\usepackage{nameref}									
\usepackage[colorlinks	=	true,						
          	 	linkcolor	=	red,
            	urlcolor	=	red,
            	citecolor	=	red]%
            	{hyperref}
\usepackage{bookmark}									
\usepackage{cleveref}									

\allowdisplaybreaks

\def\rr{{\mathbb R}}
\def\rn{{{\rr}^n}}

\def\dist{{\mathop\mathrm{\,dist\,}}}
\def\loc{{\mathop\mathrm{\,loc\,}}}

\def\R{\mathbb R}

\def\bint{{\ifinner\rlap{\bf\kern.25em--}
\int\else\rlap{\bf\kern.45em--}\int\fi}\ignorespaces}

\def\bbint{{\ifinner\rlap{\bf\kern.25em--}
\hspace{0.078cm}\int\else\rlap{\bf\kern.45em--}\int\fi}\ignorespaces}

\def\diam{{\mathop\mathrm{\,diam\,}}}

\def\r{\right}
\def\lf{\left}

\newcommand{\abs}[1]						
	{\left| #1 \right|}

\newcommand{\dtext}{\textnormal d}

\newtheorem{thm}{Theorem}[section]
\newtheorem{lem}{Lemma}[section]
\newtheorem{prop}{Proposition}[section]

\newtheorem{defn}{Definition}[section]

\newtheorem{remark}{Remark}[section]
\numberwithin{equation}{section}

\usepackage[linecolor=black!50,backgroundcolor=black!5,bordercolor=black!50,textsize=tiny]{todonotes}
\usepackage{tcolorbox}
\tcbuselibrary{breakable}

\newcommand{\onto}{\overset{{}_{\textnormal{\tiny{onto}}}}{\longrightarrow}}

\begin{document}

\arraycolsep=1pt

\title{Global Integrability of the Reciprocal of Jacobians for Homeomorphisms of Finite Distortion
}
\author{Anna Dole\v{z}alov\'{a}}
\address{Department of Mathematics and Statistics, University of
Jyv\"askyl\"a, P.O. Box 35 (MaD), 40014 Jyv\"askyl\"a, Finland
and
Department of Decision-Making Theory, Institute of Information Theory and Automation, Czech Academy of Sciences, Pod Vodárenskou věží 4, 182 00 Prague 8, Czech Republic}
\email{dolezalova@utia.cas.cz}

\author{Jani Onninen}
\address{Department of Mathematics, Syracuse University, Syracuse,
NY 13244, USA and Department of Mathematics and Statistics, University of
Jyv\"askyl\"a, P.O. Box 35 (MaD), 40014 Jyv\"askyl\"a, Finland}
\email{jkonnine@syr.edu}

\author{Yizhe Zhu}
\address{Department of Mathematics and Statistics, University of
Jyv\"askyl\"a, P.O. Box 35 (MaD), 40014 Jyv\"askyl\"a, Finland}
\email{yizhe.y.zhu@jyu.fi}

\author{Zheng Zhu}
\address{School of Mathematical science\\
         Beihang University\\
        Beijing 102206\\
        P. R. China}
\email{zhzhu@buaa.edu.cn}        



\thanks{A. Doležalová and Y. Zhu were supported by the Academy of Finland project 334014. A. Doležalová was moreover supported by the Czech Academy of Sciences project PPLZ L100752451 and by the grant GA\v CR P201/23-04766S. J.\ Onninen was supported by the NSF grant DMS-2453853. Z. Zhu was supported by the NSFC grant No. 12301111 and the Fundamental Research Funds for the Central Universities in Beihang University and the Beijing Natural Science Foundation No. 1242007.}

\date{\today}

\maketitle
\begin{abstract}
For a homeomorphism with $p$-integrable distortion, we obtain the optimal global degree of integrability for the reciprocal of its Jacobian determinant. 
As an application, we strengthen the result of  Dole\v{z}alov\'a, Hencl and Mal\'y concerning weak limits of Sobolev homeomorphisms with finite distortion. Such limits represent physically admissible deformations, as they remain injective almost everywhere and thus adhere as closely as possible to the \emph{principle of non-interpenetration of matter} in mathematical models of nonlinear elasticity.

\end{abstract}

\section{Introduction}
A homeomorphism $h \colon \Omega \to \Omega'$ of Sobolev class $W^{1,1}_{\loc} (\Omega, \R^n)$ between domains in $\R^n$ has \emph{finite distortion} if
\begin{equation}\label{eq:distineq}
\abs{Dh(x)}^n \le K(x) J_h(x) \qquad \qquad \textnormal{a.e. in } \Omega
\end{equation}
for some measurable function $1 \le K(x) < \infty$. Here $\abs{Dh(x)}$ denotes  the operator norm of the weak derivative of $h$ at $x$, and $J_h(x)=\det Dh(x)$. The smallest  $K(x)$ for which
the distortion inequality~\eqref{eq:distineq} holds, denoted by $K_h (x)$, is called the \emph{distortion function} of $h$. 

Note that the distortion inequality merely asks that the differential $Dh(x)$ vanishes whenever  the Jacobian $J_h(x)$ does, and that the Jacobian is non-negative a.e. in $\Omega$. This seems to be a minimal requirement for a
mapping to carry any  geometric information.  However, to develop a robust theory, quantitative control over the distortion is essential. In the limiting case when $K_h \in L^\infty (\Omega)$ we obtain \emph{quasiconformal} homeomorphisms.  In fact, the theory of mappings with finite distortion emerged from efforts to generalize the classical quasiconformal theory~\cite{AIMb, IMb, Reb}  in Geometric Function Theory (GFT) to the degenerate elliptic setting. This broader framework has found concrete applications  in materials science, particularly Nonlinear Elasticity (NE)~\cite{Anb, Bac, Cib} and critical phase phenomena, and in Calculus of Variations.

The general theory of hyperelasticity in the mathematical models of NE considers Sobolev homeomorphisms $h \colon \Omega \to \Omega'$
with  nonnegative  Jacobian $J_h\ge0$ that  minimize a given stored energy functional. In particular,  Neo-Hookean materials have attracted significant attention in this context (see e.g.~\cite{Ba0, BHM:2017, BHMR1, BHMR2,  BP, BPO1,  CN, conti,FG, MandS:1995, MST:1996, MTY,Sv}). The term Neo-Hookean refers to energy densities that blow up as $J_h$ tends to zero:
\begin{equation}\label{eq:neo}
\mathcal{E}^p [h] = \int_\Omega \abs{Dh(x)}^p + \phi(J_h (x))\, \dtext x,
\end{equation}
where $p \geq 1$ and $\phi$ is a positive convex function on $(0, \infty)$ with $\lim_{t\to0^+}\phi(t) = \infty $ and  $\phi(t) = \infty $ for $t \leq 0$.

It is easy to see that if $\phi(t) \gtrsim t^{-q}$, then any Sobolev homeomorphism $h \colon \Omega \to \mathbb{R}^n$ with finite energy $\mathcal{E}^p[h]$ satisfies $K_h \in L^\alpha_{\loc}(\Omega)$ for all $\alpha > 0$ such that $p > n\alpha$ and $q \geq \frac{p\alpha}{p - n\alpha}$. However, the converse is far from true. For instance, there exists a Sobolev homeomorphism $h \colon \Omega \to \Omega'$ whose Jacobian vanishes on a set of positive measure, yet it satisfies $K_h \in L^\alpha_{\loc}(\Omega)$ for all $\alpha < \frac{1}{n-1}$, see~\cite{KandM:JEMS}.  Moreover, in general, to guarantee that $J_h^{-1} \in L^q_{\loc}(\Omega)$ for some $q > 0$, one must assume that $K_h$ is essentially bounded. Indeed, for $K$-quasiconformal mappings $h \colon \Omega \to \Omega'$, Gehring’s higher integrability theorem~\cite{gehring} ensures that $J_h^{-1} \in L^q_{\loc}(\Omega)$ for some $q > 0$.
For planar $K$-quasiconformal maps, the sharp exponent $q< \frac{1}{K-1}$  was established in a seminal work by Astala~\cite{astala}. However, the sharp exponent $q=q(n,K)$ is not known when $n \ge 3$.

In this work, we address the case beyond the quasiconformal regime and establish the optimal global integrability for the reciprocal of the Jacobian determinant of a homeomorphism with $p$-integrable distortion.

\begin{thm}\label{thm:intejaco_main}
Let $\Omega\subset\rn$ be a bounded Lipschitz domain, $ \Omega'\subset\rn$ be a bounded Ahlfors $n$-regular domain, and let $h \colon \overline{\Omega}\to\overline{\Omega'}$ be a homeomorphism 
of finite distortion with  $K^{\frac{1}{n-1}}_h\in L^p(\Omega)$ for some $1\leq p<\infty$. Then there exists a constant $C$ such that we have 
\begin{equation}\label{eq:Jacointe}
\int_\Omega\log^p\lf(e+\frac{1}{J_h(x)}\r)dx\leq C\int_\Omega K^{\frac{p}{n-1}}_h(x)dx,
\end{equation}
where the constant $C=C(n, p, \Omega, \Omega', h\big|_{\partial\Omega})$ depends on $n$, $p$, $\Omega$, $\Omega'$ and also the boundary value of $h$ on $\partial\Omega$.
\end{thm}

The local version of Theorem~\ref{thm:intejaco_main} was established by Koskela, Onninen, and Rajala in~\cite{KOR:JGA}, where they also provided an example showing that the degree of integrability of the distortion function in~\eqref{eq:Jacointe} cannot be lowered without losing the conclusion; see~\cite[Example 1.1]{KOR:JGA}.

Integral bounds of this type naturally arise, for instance, in connection with refined forms of Morrey’s quasiconvexity condition in two dimensions, under which one can prove lower semicontinuity and the existence of minimizers; see~\cite{AFGKK, AFGKK2}.

The Lipschitz regularity assumption on the reference configuration $\Omega$  is a standard, axiomatic condition in NE. On the other hand, we assume the target domain $\Omega '$ only to be a bounded Ahlfors $n$-regular domain. Such domains, frequently encountered in GFT, admit irregular, possibly fractal boundaries, yet retain crucial geometric measure properties such as uniform boundary measure scaling~\cite{He_metric, HKActa, HKST}. This broader framework still allows one to pass from local to global integrability results, thereby encompassing a much wider class of geometric configurations than the smooth or Lipschitz category.

Nevertheless, one cannot expect the global estimate~\eqref{eq:Jacointe} to hold with a constant $C$ independent of $h$ even for smooth target domains. Indeed, consider the  sequence of M\"obius transformations $\{h_k\}_{k=1}^\infty$ from the unit disk $\mathbb D\subset \mathbb{C}$ onto itself given by the formula 
\[h_k(z)=\frac{z-a_k}{a_k z-1} \, , \qquad \qquad\textnormal{ where }\quad  a_k=1-\frac{1}{2k}\, .\]
Then
\[\lim_{k\to \infty} \int_{\mathbb D} \phi (J_{h_k}(z)) \, dz = \infty\]
for any $\phi$ which is a nonnegative decreasing  function on $(0, \infty)$ with $\lim_{t\to0^+}\phi(t)=\infty$. However, as $h_k$ is conformal, $K_{h_k}\equiv 1$ for each $k$. In particular, 
$$\int_{\mathbb D}K_{h_k}^{\frac{p}{n-1}}(z)\, dz =\pi .$$

To study minimization problems, one must first address the question which class of admissible mappings to consider. Even for the classical Dirichlet energy in two dimensions, injectivity may be lost in the limit of minimizing sequences of homeomorphisms, leading to \emph{interpenetration of matter}~\cite{IKKO, IOlag,IandO:JEMS, IOhopf}. Therefore, the development of viable  models in both 2D and 3D demands a systematic study of weak limits of Sobolev homeomorphisms and, more generally, a principled framework for enlarging the class of admissible deformations.

We denote by $\mathscr H^p_\varphi  (\overline{\Omega}, \overline{\Omega'})$ the class of homeomorphisms $h \colon \overline{\Omega} \onto \overline{\Omega'}$ in the Sobolev space $W^{1,p} (\Omega, \mathbb R^n)$ that coincide with a given boundary homeomorphism $\varphi \colon \partial \Omega \onto \partial \Omega'$ on $\partial \Omega$. We assume that this class is nonempty.
 For details on conditions under which $\varphi$ admits a Sobolev homeomorphic extension (i.e., when $\mathscr H^p_\varphi (\overline{\Omega}, \overline{\Omega'}) \neq \emptyset$), we refer the reader to~\cite{HKOext, KKO,KO1, KOsjsp}.

If $p\ge n$, then a sequence of homeomorphisms $h_k \in \mathscr H_\varphi^p (\overline{\Omega}, \overline{\Omega'})$ contains an uniformly convergent subsequence, and therefore the limit mapping
is monotone. \emph{Monotonicity} of a continuous map $h \colon \overline{\Omega} \onto \overline{\Omega'}$, the concept of  Morrey~\cite{Mor}, simply means that  the preimage $h^{-1} (y_0)$ of any point  $ y_0 \in \overline{\Omega'}$ is a continuum in $\overline{\Omega}$.  In the planar case the class of $W^{1,p}$-monotone mappings $h \colon  \overline{\Omega} \onto \overline{\Omega'}$
 which coincide with $\varphi$ on $\partial \Omega$  actually fully characterizes  the class of $W^{1,p}$-weak limits of homeomorphisms when $p\ge 2$.
 As a consequence, weak and strong limits of planar $W^{1,p}$-homeomorphisms coincide when $p \geq 2 = n$, see also~\cite{IOmono}. A deep result of De Philippis and Pratelli~\cite{DP} extends   this to the limits of self-homeomorphisms between closed unit disks that equal the identity mapping on the boundary when $p<2$. 
 
The situation, however, changes drastically when $p < n$. In this case, weak limits of  $W^{1,p}$-homeomorphisms need not be continuous, allowing for phenomena such as cavitation. Consequently, monotone maps no longer suffice as admissible deformations. In models of NE, however, only physically unacceptable maps are typically excluded, most importantly, those exhibiting \emph{strong interpenetration of matter}, or where volumes of material turn ``inside out''.  Motivated by this, M\"uller and Spector~\cite{MandS:1995}, building on ideas of \v Sver\'ak~\cite{Sv}, introduced the notion of Sobolev $(INV)$ mappings. Informally, an  $(INV)$ map sends a ball $B(x, r)$ inside its own boundary image $h(\partial B(x, r))$, while mapping the complement of the closed ball outside this boundary; see Section~\ref{sect_prelim} for a precise formulation. This condition was originally defined for $W^{1,p}$-mappings with $p > n-1$, where continuity on almost every sphere ensures the notion of "inside" and "outside" is well-defined. For $p > n-1$, weak limits of homeomorphisms in $\mathscr H^p_\varphi (\overline{\Omega}, \overline{\Omega'})$ satisfy the $(INV)$ condition, and moreover, the class of such $(INV)$ maps is weakly closed in $W^{1,p}$~\cite[Lemma 3.3]{MandS:1995}. If such deformations also have strictly positive Jacobian determinant a.e.  (i.e., are \emph{strictly orientation-preserving}), they are injective a.e.~\cite{MandS:1995}. 
 
There has been a growing interest in extending these results to the borderline case of $p=n-1$.  The inception occurred when Conti and De Lellis~\cite{conti} were able to introduce the $(INV)$ condition also  for bounded mappings in the Sobolev class $W^{1, n-1}$. Their studies were mainly motivated by the 3D models in  compressible Neo-Hookean elasticity which contains the classical Dirichlet term $|Dh|^2$, i.e., mappings $h$ with $\mathcal E^2 [h] < \infty$ (see~\eqref{eq:neo}). Strictly orientation-preserving bounded $W^{1,n-1}$ mappings that satisfy the $(INV)$ condition  share many important properties (such as being injective a.e.) with  the corresponding mappings in  $W^{1,p} (\Omega, \R^n)$ when $p>n-1$.  However, when $p=n-1$, this class is not closed under weak limits. Conti and De Lellis~\cite{conti} constructed a sequence of bi-Lipschitz mappings on the unit ball $\mathbb{B} \subset \mathbb{R}^3$, which coincide with the identity on the boundary and have uniformly bounded energy $\mathcal E^2$, which converges weakly to a mapping violating the $(INV)$ condition. Doležalová, Hencl, and Malý~\cite{INV} further showed that such a limit cannot be a strong limit of $W^{1,2}$-homeomorphisms, thereby highlighting a fundamental difference from the planar case: in 3D, weak and strong limits of $W^{1,2}$-homeomorphisms do not coincide.

In~\cite{INV}, conditions under which weak $W^{1,n-1}$-limits of homeomorphisms satisfy the $(INV)$ condition were identified. Specifically, suppose that the function $\phi $ in \eqref{eq:neo} satisfies moreover that there exists $A > 0$ such that
$$
A^{-1} \phi(t) \leq \phi(2t) \leq A \phi(t) \quad \text{for all } t \in (0, \infty).
$$
Then if a sequence of $W^{1, n-1}$-homeomorphisms ${h_k}$ between closures of bounded Lipschitz domains shares the same boundary data and satisfies
\begin{itemize}
\item[(A)] $\sup_k \mathcal E^{n-1} [h_k] < \infty $, and
\item[(B)] $\sup_k \int_{\Omega} K^\frac{1}{n-1}_{h_k} < \infty$,
\end{itemize}
 then its weak limit satisfies the $(INV)$ condition. 
 
By using Theorem~\ref{thm:intejaco_main}, we show that for sequences weakly converging in $W^{1, n-1}$, the assumption (A) is redundant as they are bounded and we can estimate the integral containing $\phi$ by the distortion term:
\begin{thm}\label{thm:outdistortion}
Let $\Omega, \Omega'\subset\rn$ be bounded Lipschitz domains and $n \ge 3$. Suppose that $h_k \colon \overline{\Omega} \onto \overline{\Omega'}$ are homeomorphisms in $W^{1, n-1}(\Omega, \rn)$ with $h_k\equiv h_1$ on $\partial\Omega$ for all $k\in\mathbb N$ and
\begin{equation}\label{eq:supremum}
\sup_{k}\int_{\Omega}K_{h_k}^{\frac{1}{n-1}}(x)dx<\infty.
\end{equation}
If the sequence $\{h_k\}$ converges weakly in $W^{1, n-1}(\Omega, \rn)$   to a mapping $h$, then  $h$ satisfies $(INV)$.
\end{thm}
The  $L^\frac{1}{n-1}$-integrability assumption for the distortion functions in (\ref{eq:supremum}) is sharp, see Remark~\ref{rk:integrabilty_K}.

\section{Preliminary}\label{sect_prelim}
In this section, we give definitions and results which are needed in the following sections. In what follows, $B(z,r)$ denotes an open ball with center $z\in\rn$ and radius $r>0$. If we do not need to emphasize the center, we also use $B_r$ to denote a ball with the radius $r>0$. We write $C$ as a generic positive constant, which may change even in a single string of estimates. The dependence of a constant on parameters $\alpha, \beta,\cdots$ is expressed by the notation $C = C(\alpha, \beta, \cdots)$ if needed to be emphasised. We start by recalling the definition of Lipschitz and Ahlfors $n$-regular domains, for more see e.g.~\cite{Evans}.

\subsection{Lipschitz and Ahlfors regular domains}

\begin{defn}\label{de:Lipschitz}
A bounded domain $\Omega\subset\rn$ is called {\it Lipschitz}, if for every $z\in\partial\Omega$ there exists $r_z>0$ and a Lipschitz mapping $M_z:\rr^{n-1}\to\rr$ such that, upon rotating and relabeling the coordinate axes if necessary, we have
\[\Omega\cap B(z, r_z)=\lf\{x=(x_1, x_2, \cdots, x_n):M_z(x_1, x_2, \cdots, x_{n-1})<x_n\r\}\cap B(z, r_z).\]
\end{defn}
In other words, near every point $z\in\partial\Omega$ the boundary of $\Omega$ is a graph of a Lipschitz function. For every $r>0$, we define the upper half-ball by 
\[B^+(0, r):=\lf\{x=(x_1, x_2,\dots, x_n)\in\rn: |x|<r\ {\rm and}\ x_n>0\r\}\]
and the projection of the ball $B(0, r)$ to the $(n-1)$-dimensional hyperplane $\rr^{n-1}\times\{0\}$ by 
\[\mathcal{P}B(0, r):=\lf\{x=(x_1,x_2,\dots, x_n)\in\rn:|x|<r\ {\rm and}\ x_n=0\r\}.\]
The following property of Lipschitz domains is well-known, see e.g. \cite{Licht, monk}. 

\begin{lem}\label{le:bilipschitz}
Let $\Omega\subset\rn$ be a Lipschitz domain. Then for every $z\in\partial\Omega$ there exist $r_z>0$ and a bi-Lipschitz homeomorphism $L_z \colon B(z, r_z) \onto B(0, r_z)$ such that
$$L_z: {\Omega\cap B(z, r_z)}\onto {B^+(0, r_z)} \, , $$
$L_z(z)=0$, and 
\[L_z\lf(\partial\Omega\cap B(z, r_z)\r)=\mathcal{P}B(0, r_z).\]
\end{lem}

Now, let us define Ahlfors $n$-regular domains.
\begin{defn}\label{de:regular}
A domain $\Omega\subset\rn$ is called an Ahlfors $n$-regular domain, if there exists a constant $0<c<1$ such that for every $x\in\overline\Omega$ and $0<r<1$ we have 
\[\mathcal H^n(\Omega\cap B(x ,r))\geq c\mathcal H^n(B(x, r)).\]
\end{defn}

\begin{remark}
Given a bounded Lipschitz domain $\Omega$, $x\in \overline{\Omega}$ and $r\in (0,\diam(\Omega))$, we have that
\begin{equation}\label{lipschitz_ahlfors}
\frac{r^n}{C}\leq |B(x,r)\cap \Omega|\leq C r^n
\end{equation}
where the constant $C>0$ depends on the shape of $\Omega$. This means that any bounded Lipschitz domain is an Ahlfors $n$-regular domain.
\end{remark}

\subsection{The change-of-variables formula}

The following change-of-variables formula is a special case of \cite[Theorem 3.1.8]{Federer}, or~\cite[Theorem A.35]{HandK:book}.
\begin{lem}\label{le:area}
  Let $\Omega,\Omega'\subset\rn$ be domains, and let $h:\Omega\to\Omega'$ be a homeomorphism  in the Sobolev  class $W_{\loc}^{1,1}(\Omega, \Omega')$.  Then for every nonnegative Borel measurable function $\eta$ on $\Omega'$ we have  
  \begin{equation}\label{eq:area_0}
 \int_{\Omega}\eta(h(x))|J_h(x)|dx\leq \int_{\Omega'}\eta(y) dy.
 \end{equation}
Furthermore, there exists a full-measure  subset $\widetilde\Omega\subset \Omega$ such that we have an identity in ~\eqref{eq:area_0}; that is, 
  \begin{equation}\label{eq:area}
 \int_{\widetilde\Omega}\eta(h(x))|J_h(x)|dx=\int_{h(\widetilde\Omega)}\eta(y)dy.
 \end{equation}

\end{lem}

\subsection{The co-area formula}
The following co-area formula for Sobolev functions is well-known, see e.g.~\cite[Theorem 1.1]{MSZ:tams}. 
\begin{prop}\label{prop:co-area}
Let  $u \colon \Omega \to \R$ be a function in $W^{1,1} (\Omega, \mathbb{R} )$. Then for every measurable set $E \subset \Omega$ we have
\[
\int_E \abs{ \nabla  u (x) } \, dx = \int_{- \infty}^{\infty} \mathcal H^{n-1}(E \cap u^{-1} (s)) \, ds.
\]
\end{prop}
Here $\mathcal H^{n-1} (A)$ denotes the $(n-1)$-dimensional Hausdorff measure of the set $A$.

\subsection{The $(INV)$ condition}  
In this section, we introduce the $(INV)$ condition for the convenience of the reader. For a more detailed and formal presentation see \cite{conti} or \cite[Section 2.4]{INV}.

Let $\Omega\subset\rn$ be a bounded domain and $f:\overline\Omega\to\rn$ be a smooth mapping. Then the (classical) topological degree of $f$ at $y$ with respect to $\Omega$ is defined as
\[{\rm deg}(f, \Omega, y):=\sum_{x\in\Omega\cap f^{-1}(y)}{\rm sgn}(J_f(x))\]
if $J_f(x)\neq 0$ for each $x\in\Omega\cap f^{-1}(y)$. With an approximation argument, this definition can be extended to an arbitrary continuous mapping from $\overline\Omega$ to $\rn$. We should emphasize that the degree only depends on values of $f$ on $\partial\Omega$, see \cite{degree:book} for more details. For a homeomorphism $f :\overline\Omega\to\rn$, either ${\rm deg} (f, \Omega, y)=1$ for all $y\in f(\Omega)$ ($f$ is sense-preserving), or ${\rm deg}(f, \Omega, y)=-1$ for all $y\in f(\Omega)$ ($f$ is sense-reserving). For a $W^{1, n-1}$-homeomorphism, the sign of the Jacobian shows the topological orientation, see the following proposition from \cite{HandM:CVPDE}.

\begin{prop}\label{prop:degree}
  Let $f\in W^{1, n-1}(\Omega, \rn)$ be a homeomorphism on $\overline\Omega$ with $J_f>0$ almost everywhere. Then 
  \[{\rm deg}(f, \Omega, y)=1\ \ {\rm for\ every}\ y\in f(\Omega).\]
\end{prop}

The notion of the topological degree allows us to talk about "inside" of a continuous image of a sphere by taking the points where the degree is non-zero. For mappings in $W^{1,p}(\Omega, \mathbb{R} )$, $p>n-1$, we know that they are continuous on the sphere $S(a,r)$ for every $a\in \Omega$ and a.e. $r\in(0, \dist(a,\partial \Omega))$. Therefore, we can use the classical degree. However, in the borderline case $p=n-1$, one needs a generalized version of the degree, as our functions may exhibit discontinuities on a lot of spheres. The function $\rm Deg$ is the (unique) integer-valued function in $BV(\mathbb{R}^n)$ such that 
\begin{equation}\label{eq:Degree}
    \int_{\rn}{\rm Deg}(f, B, y)\psi(y)dy=\int_{\partial B}(\mathbf{u}\circ f)\cdot(\Lambda_{n-1} D_\tau f)\mathbf{v}d\mathcal H^{n-1}
\end{equation}
for every test function $\psi\in C_0^\infty(\rn)$ and every $C^\infty$-vector field $\mathbf u$ on $\rn$ satisfying ${\rm div}[\mathbf u]=\psi$,
where $\mathbf{v}(x)$ denotes the exterior normal vector to $B$ at $x\in\partial B$ and $D_\tau f$ the tangential derivative of $f$ on $\partial B$. 
One can prove that for $f\in W^{1, n-1}(\partial B, \rn)\cap C(\overline B, \rn)$ with $|f(\partial B)|=0$, we have
\[{\rm Deg}(f, B, y)={\rm deg}(f, B, y)\ {\rm for\ almost\ every\ }y\in\rn.\]

\begin{defn}\label{de:Topimage}
 Given a measurable set $E\subset\rn$, we call 
\[\lim_{r\to0^+}\frac{\mathcal H^n(E\cap B(x, r))}{\mathcal H^n(B(x, r))}\]
the density of the set $E$ at the point $x\in\rn$, if the limit exists.

Let $B\subset\rn$ be a ball and $f\in W^{1, n-1}(\partial B, \rn)\cap L^\infty(\partial B, \rn)$. We define $im_T(f, B)$, the topological image of $B$ under $f$, as the set of all points where the density of the set $\{y\in\rn:{\rm Deg}(f, B, y)\neq 0\}$ is $1$.
\end{defn}
Let $\Omega\subset\rn$ be a domain. By the notation $U\subset\subset\Omega$ we mean that $\overline U\subset\Omega$. Now, we are ready to define the $(INV)$ condition as in \cite{conti, MandS:1995}.
\begin{defn}\label{de:INV}
Let $f\in W^{1, n-1}(\Omega, \rn)\cap L^\infty(\Omega, \rn)$. We say that $f$ satisfies the $(INV)$ condition in the ball $B\subset\subset\Omega$ if
\begin{itemize}
\item[$(1)$] its trace on $\partial B$ is in $W^{1, n-1}\cap L^\infty$;

\item[$(2)$] $f(x)\in im_T(f, B)$ for almost every $x\in B$;

\item[$(3)$] $f(x)\notin im_T(f, B)$ for almost every $x\in\Omega\setminus B$.\\
\end{itemize}
We say that $f$ satisfies the $(INV)$ condition in $\Omega$ if for every $a\in\Omega$ there is $r_a>0$ such that for $\mathcal H^1$-a.e. $r\in(0, r_a)$ it satisfies the $(INV)$ condition in $B(a, r)$.
\end{defn}

\section{Proof of Theorem \ref{thm:intejaco_main}}
In this section, we prove  Theorem~\ref{thm:intejaco_main}.

We first fix $w\in\partial\Omega'$ and set $z:=h^{-1}(w)$.  Then we  choose sufficiently small radii $r_w>0$ and $t_z>0$  so that the bi-Lipschitz mapping from Lemma~\ref{le:bilipschitz} exists, i.e.,
$$ 
L_z:\overline{\Omega\cap B(z, t_z)}\to\overline{B^+(0, t_z)},
$$ 
and
that
$$
h^{-1}(B(w, r_w)\cap \partial \Omega')\subset B(z, t_z)\cap \partial \Omega.
$$
Note that the required smallness of the radii depends on the domains $\Omega$ and $\Omega'$, and the boundary map $h \colon \partial \Omega \to \partial \Omega'$.

 \begin{lem}
For every ball $B_r$ with $B_r\cap\Omega'\subset B(w, r_w)\cap\Omega'$, we have 
\begin{equation}\label{eq:s-iso}
\mathcal H^{n-1}\lf(h^{-1}\lf(\partial\Omega'\cap B_r\r)\r)\leq C \mathcal H^{n-1}\lf(h^{-1}\lf(\partial B_r\cap\Omega'\r)\r),
\end{equation}
where the constant C depends only on the bi-Lipschitz constant of the map $L_z$.
\end{lem}
\begin{proof}
For every $B_r$ with $B_r\cap\Omega'\subset B(w,r_w)\cap\Omega'$, the set $h^{-1}(B_r\cap\partial\Omega')$ is contained in $\overline{B(z,t_z)\cap\Omega}$, which is $L_z$-bi-Lipchitz equivalent to $\overline{B^+(0,t_z)}$. In particular, the image $A:=L_z(h^{-1}(B_r\cap\partial\Omega'))$ is contained in the flat base of $\overline{B^+(0,t_z)}$, and
\begin{align}\label{s-iso-1}
\mathcal{H}^{n-1}(h^{-1}(B_r\cap\partial\Omega'))\le C(L_z)\mathcal{H}^{n-1}(A).
\end{align}
Let $A':=L_z(h^{-1}(\partial B_r\cap\Omega'))$ be the image surface and denote $\mathcal{P}A'$ by its projection on the base of $\overline{B^+(0,t_z)}$. Note that
$A\subset \mathcal{P}A'$ and
\begin{align*}
\mathcal{H}^{n-1}(A')=\int_{\mathcal{P}A'}\sqrt{\left(\frac{\partial L_z}{\partial x_1}\right)^2+\cdots+\left(\frac{\partial L_z}{\partial x_{n-1}}\right)^2+1}\ d\mathcal{H}^{n-1}.
\end{align*}
Therefore,
\begin{align}\label{s-iso-2}
    \mathcal{H}^{n-1}(A)\le\mathcal{H}^{n-1}(\mathcal{P}A')\le \mathcal{H}^{n-1}(A')\le C(L_z)\mathcal{H}^{n-1}(h^{-1}(\partial B_r\cap\Omega')).
\end{align}
Combining \eqref{s-iso-1} and \eqref{s-iso-2}, we get
\begin{align}\label{s-iso}
   \mathcal{H}^{n-1}(h^{-1}( B_r\cap\partial\Omega'))\le C(L_z)\mathcal{H}^{n-1}(h^{-1}(\partial B_r\cap\Omega')). 
\end{align}
\end{proof}

 Since $\partial\Omega'$ is compact, there exists an integer $M_1\in\mathbb N$ such that 
\begin{equation}\label{eq:w1}
\partial\Omega'\subset\bigcup_{i=1}^{M_1}B\lf(w_i, \frac{r_{w_i}}{2}\r) = \bigcup_{i=1}^{M_1}\frac{1}{2}B\lf(w_i, r_{w_i}\r).
\end{equation}
Since $\overline{\Omega'}\subset\rn$ is also a compact set, there exists a finite collection of balls $\{B^k\subset\Omega'\}_{k=1}^{M_2}$ such that
\begin{equation}\label{eq:w2}
 \overline{\Omega'} \subset\bigcup_{i=1}^{M_1}\frac{1}{2}B\lf(w_i, r_{w_i}\r)\cup\bigcup_{k=1}^{M_2}\frac{1}{2}B^k
\end{equation}
and 
\begin{equation}\label{eq:w3}
\sum_{i=1}^{M_1}\chi_{B(w_i, r_{w_i})}+\sum_{k=1}^{M_2}\chi_{B^k}\leq C(n)<\infty \, . 
\end{equation}
We denote such a family of balls by $\mathcal B=\mathcal B_1 \cup \mathcal B_2$, where 
\[ \mathcal B_1 = \{ B(w_i, r_{w_i}) \colon i=1, \dots , M_1\}  \textnormal{ and }  \mathcal B_2 = \{ B^k \colon k=1, \dots , M_2\} . \]  
Note that the choice of $\mathcal B$ is universal, given $\Omega$, $\Omega'$ and the boundary values of $h$.

Let $h:\overline{\Omega}\to\overline{\Omega'}$ be a homeomorphism such that its restriction to $\Omega$ lies in $W^{1,1}_{\rm loc}(\Omega, \Omega')$ and has finite distortion with $K_h^{\frac{p}{n-1}}\in L^1(\Omega)$ for some $1\leq p<\infty$. Let $N$ be a set of measure zero such that  $h$ satisfies the change-of-variables formula~\eqref{eq:area} on $\Omega \setminus N$. We define two functions $j$ and $g$ on $\Omega'$,  as in \cite{KOR:JGA},  by setting 
\begin{equation}\label{eq:j}
    j(y):=\begin{cases}
    \frac{1}{J_h(h^{-1}(y))},\ \ & {\rm if}\ J_{h}(h^{-1}(y))>0\ {\rm and}\ h^{-1}(y)\in \Omega \setminus N,\\
    0,\ \ & {\rm elsewhere},
    \end{cases}
\end{equation}
and 
\begin{equation}\label{eq:psi}
g(y):=\begin{cases}\frac{|Dh(h^{-1}(y))|}{J_h(h^{-1}(y))},\ \ & {\rm if}\ J_h(h^{-1}(y))>0\ {\rm and}\ h^{-1}(y)\in \Omega \setminus N, \\
    0,\ \ & {\rm elsewhere}.
    \end{cases}
\end{equation} 
By  \cite[Theorem 1.2]{KandM:JEMS}, for almost every $x\in\Omega$ we have $J_h(x)>0$. We denote
\[
\widetilde \Omega' :=\lf\{y\in \Omega' \colon  j(y)>0\r\}  \, . 
\]
Note that if $h$ is differentiable a.e., then for a.e. $y\in \widetilde \Omega'$ we have that $j(y)=J_{h^{-1}}(y)$ and $g(y)=|D^\#h^{-1}(y)|$, where $D^\#$ denotes the cofactor matrix. By (\ref{eq:j}), we have $h^{-1}(\widetilde \Omega')\subset \Omega \setminus N$. Let $U$ be a measurable subset of $\widetilde \Omega'$.  Hence, we can apply (\ref{eq:area}) to $h^{-1}(U)$ to obtain
\begin{equation}\label{eq:area2}
\int_{h^{-1}(U ) }\eta(h(x))J_h(x)dx=\int_{U }\eta(y)dy
\end{equation}
for every nonnegative measurable function $\eta$ defined on $\Omega'$. 

A crucial estimate that enabled us to prove the local version of Theorem~\ref{thm:intejaco_main} (in~\cite{KOR:JGA}) was  the following integral variant of the inverse isoperimetric inequality 
\begin{equation}\label{eq:NEW1}
\bint_{B_r}j(y)dy\leq C(n)\lf(\bint_{B_{2r}}g(y)dy\r)^{\frac{n}{n-1}},
\end{equation}
where $B_r$ is an arbitrary ball such that $B_{2r} \subset \Omega'$.  Throughout this text, we use the notation $\bint_A = \frac{1}{|A|} \int_A$ to denote the integral average over the set $A$.

The following lemma extends this formula to include sufficiently small balls that intersect the boundary of $\Omega'$. 
\begin{lem}\label{le:IIE}
Let $B \in \mathcal B_1$ and let  $B_r\subset\rn$ be a ball whose center $y_0$ is contained in $B\cap\Omega'$ with ${B_{2r}}\cap\Omega'\subset B \cap\Omega'$.  Then 
\begin{equation}\label{eq:IIE}
    \bint_{B_r\cap\Omega'}j(y)dy\leq C\lf(\bint_{B_{2r}\cap\Omega'}g(y)dy\r)^{\frac{n}{n-1}},
\end{equation}
where the constant $C$ depends on $n$, $\Omega$ and $\Omega'$.
\end{lem}
\begin{proof}
Clearly, 
\begin{equation}\label{equa1}
 \int_{B_r\cap\Omega'}j(y)dy= \int_{\{y\in B_r\cap\Omega' \colon j(y)=0 \}}j(y) dy + \int_{{B_r\cap \widetilde \Omega'}}j(y)dy =  \int_{B_r\cap \widetilde \Omega'}j(y)dy \, .
 \end{equation}
  By (\ref{eq:area2}), we have
 \begin{equation}\label{equa2}
 \begin{split}
 \int_{{B_r\cap \widetilde \Omega'}}j(y)dy &=\int_{h^{-1}\lf({B_r\cap \widetilde \Omega'}\r)}j(h(x))J_h(x)dx
=\int_{h^{-1}\lf({B_r\cap \widetilde \Omega'}\r)}\frac{1}{J_h(x)}J_h(x)dx\\
& = \abs{h^{-1}\lf({B_r\cap \widetilde \Omega'}\r)} \leq  \abs{ h^{-1}(B_r\cap\Omega')} .
\end{split}
\end{equation}
 For every $s\in(r, 2r)$, we distinguish two cases. If $B_s\subset B \cap\Omega'$, the classical isoperimetric inequality gives 
\[
     \abs{ h^{-1}(B_s)}^{\frac{n-1}{n}}\leq C(n)\,  \mathcal H^{n-1}\lf( \partial ( h^{-1}\lf( B_s\r) ) \r) 
     = C(n)\,  \mathcal H^{n-1}\lf(  h^{-1}\lf(\partial B_s\r)  \r)  ,  \nonumber 
\]
because $h$ is a homeomorphism. 
If   $B_s\cap\partial\Omega'\neq\emptyset$, the classical isoperimetric inequality and the inequality (\ref{eq:s-iso}) imply
\[
\begin{split}
\abs{ h^{-1}(B_s \cap {\Omega'}) }^{\frac{n-1}{n}} &\leq C \, \mathcal H^{n-1}\lf(\partial(h^{-1}( B_s\cap{\Omega'}))\r)\nonumber\\
&\leq C \, \mathcal H^{n-1}\lf(h^{-1}\lf(\lf(\partial\Omega'\cap B_s\r)\cup\lf(\partial B_s\cap\Omega'\r)\r)\r)\nonumber\\
&\leq C \,  \mathcal H^{n-1}\lf(h^{-1}\lf(\partial B_s\cap\Omega'\r)\r),
 \end{split}
 \]
where the constant $C$ depends on $n$,  $\Omega$ and $\Omega'$, but it doesn't depend on $B_s$.
Therefore, we always have the following isoperimetric-type inequality
\[
 \abs{ h^{-1}(B_s\cap {\Omega'})}^{\frac{n-1}{n}}  \leq C(n, \Omega, \Omega') \, \mathcal H^{n-1}\lf(h^{-1}\lf(\partial B_s\cap\Omega'\r)\r).
\]
By integrating this inequality over the interval $(r, 2r)$ we obtain
  \begin{equation}\label{eq:joo} 
  \begin{split}
  r  \abs{ h^{-1}(B_r\cap {\Omega'})}^{\frac{n-1}{n}} &\leq \int_r^{2r}\abs{ h^{-1}(B_s\cap {\Omega'})}^{\frac{n-1}{n}} \, ds  \\
   &\leq C \,\int_r^{2r} \mathcal H^{n-1}\lf(h^{-1}\lf(\partial B_s\cap\Omega'\r)\r)  \, ds.
   \end{split}
 \end{equation}
 On the other hand, applying the co-area formula (Proposition~\ref{prop:co-area}) to $u=|h-y_0|$ and $E=h^{-1}(B_{2r}\cap \Omega')$ gives
\[
\begin{split}
 \int_r^{2r} \mathcal H^{n-1}\lf(h^{-1}\lf(\partial B_s\cap\Omega'\r)\r)  \, ds &= \int_r^{2r} \mathcal H^{n-1}\lf(E\cap h^{-1}\lf(\partial B_s\cap\Omega'\r)\r)  \, ds \\
 &\leq  \int_{h^{-1}(B_{2r}\cap\Omega')}  \abs{ \nabla u(x)} \, dx
\leq \int_{h^{-1}(B_{2r}\cap\Omega')} \abs{Dh(x)} \, dx.                
\end{split} 
 \]
 Combining this with~\eqref{eq:joo}, we have
 \begin{equation}\label{eq:joo_too}
    r  \abs{ h^{-1}(B_r\cap \Omega')}^{\frac{n-1}{n}} \leq  C(n, \Omega, \Omega') \,\int_{h^{-1}(B_{2r}\cap\Omega')} \abs{Dh(x)} \, dx.
    \end{equation}
 Since $h$ is a homeomorphism with $K_h \in L^{\frac{1}{n-1}} (\Omega)$, it follows from  \cite{KandM:JEMS} that the Jacobian of $h$ is positive almost everywhere in $\Omega$, and therefore
 \[
\int_{h^{-1}\lf(B_{2r}\cap\Omega'\r)}|Dh(x)|dx=\int_{h^{-1}(B_{2r}\cap\Omega')}\frac{|Dh(x)|}{J_h(x)}J_h(x)dx\\
\leq\int_{h^{-1}(B_{2r}\cap\Omega')}g(h(x))J_h(x)dx.
\]
 Further, the change-of-variables formula~\eqref{eq:area_0} gives
\begin{equation}\label{equa4}
 \int_{h^{-1}\lf(B_{2r}\cap\Omega'\r)}|Dh(x)|dx\leq\int_{B_{2r}\cap\Omega'}g(y) \, dy.
  \end{equation}
From \eqref{lipschitz_ahlfors}, we know that
 $$
 \frac{r^n}{C} \leq |B_{r}\cap\Omega' |\leq |B_{2r}\cap\Omega' |\leq C  r^n.
 $$
  By combining this with inequalities \eqref{equa1}, (\ref{equa2}), (\ref{eq:joo_too}) and (\ref{equa4}),
   we obtain the claimed inequality~\eqref{eq:IIE}.
\end{proof}

The basic idea of the following proof comes from~\cite{KOR:JGA}.
\begin{lem}\label{le:KRO}
Under the assumptions of Theorem \ref{thm:intejaco_main}, for every $B\in\mathcal{B}$ we have 
\begin{equation}\label{eq:ball}
\int_{h^{-1}\lf(\frac{1}{2}B\cap\Omega'\r)}\log^p\lf(e+\frac{1}{J_h(x)}\r)dx\leq C\int_{h^{-1}\lf(B\cap\Omega'\r)}K_h^{\frac{p}{n-1}}(x)dx, 
\end{equation}
where the constant $C>0$ depends on $n, p$, $|\Omega|$, $|B\cap\Omega'|$, the boundary value of $h$ on $\partial\Omega$ and the constant from \eqref{lipschitz_ahlfors}.
\end{lem}
\begin{proof}
The key estimate for this proof is
\begin{equation}\label{eq:key_ineq}
\int_{\frac{1}{2}B\cap\Omega'}j(y)\log^p\lf(e+\frac{j(y)}{j_{B\cap\Omega'}}\r)dy
\leq C\int_{B\cap\widetilde\Omega'}\lf[g(y)\r]^{\frac{np}{n-1}}\lf[j(y)\r]^{1-p}dy,
\end{equation} 
where $j_{B\cap\Omega'}$ is the integral average of $j$ over the set $B\cap\Omega'$.
For $B\in\mathcal{B}_2$, this was proved in Lemma 2.4 of~\cite{KOR:JGA}, where the constant $C$ depends only on $n$ and $p$. For $B\in \mathcal{B}_1$, we provide a proof below (with more dependencies). Note that the inequality \eqref{eq:key_ineq} is invariant with respect to rescaling of $h$. Indeed, assume $S>0$ and let $\bar{h}:S\Omega\to\Omega'$ be defined as $\bar{h}(x):=h(S^{-1}x)$. For such function we have $\bar{j}(y)=S^{n}j(y)$ and $\bar{g}(y)=S^{n-1}g(y)$. That allows us to assume 
\begin{equation}\label{eq:normalization}
\widetilde{C}\left(\int_{B\cap\Omega'}j(y)dy\right)^{\frac{1}{n}}=\frac{1}{2},
\end{equation}
since we can set $S$ appropriately and arrive to \eqref{eq:key_ineq} for $\bar{h}$, which then implies it also for the original function $h$.

We first show that once we have this estimate, the statement of the lemma follows: 
First combine the definition of the function $j$ with (\ref{eq:area}), which yields
\begin{equation*}
\begin{split}
\int_{B\cap\Omega'}j(y)dy& =\int_{B\cap \widetilde \Omega'}j(y)dy=\int_{h^{-1}(B\cap\widetilde \Omega')}j(h(x))J_h(x)dx
\\ & =\int_{h^{-1}(B\cap\widetilde \Omega')}\frac{1}{J_h(x)}J_h(x)dx
=\abs{ h^{-1}(B\cap\widetilde \Omega')} \leq\abs{ \Omega }.
\end{split}
\end{equation*}
This gives an auxiliary inequality
\begin{align*}
\int_{h^{-1}\lf(\frac{1}{2}B\cap\Omega'\r)}\log^p\lf(e+\frac{1}{J_h(x)}\r)&\leq C\int_{h^{-1}\lf(\frac{1}{2}B\cap\Omega'\r)}\log^p\lf(e+\frac{\abs{B\cap\Omega'}}{\abs{\Omega}J_h(x)}\r)dx \\
&\leq  C\int_{h^{-1}\lf(\frac{1}{2}B\cap\Omega'\r)}\log^p\lf(e+\frac{1}{j_{B\cap\Omega'} \cdot J_h(x)}\r)dx,
\end{align*}
where the constant $C$ depends only on $p$ and the ratio $|B\cap\Omega'|/|\Omega|$ (which is bounded from below as the collection $\mathcal{B}$ is finite).
Recall that the Jacobian of $h$ is positive a.e., and therefore the change-of-variables formula (\ref{eq:area_0}) gives 
\begin{equation*}
\int_{h^{-1}\lf(\frac{1}{2}B\cap\Omega'\r)}\log^p\lf(e+\frac{1}{j_{B\cap\Omega'} \cdot J_h(x)}\r)dx\leq \int_{\frac{1}{2}B\cap\Omega'}j(y)\log^p\lf(e+\frac{j(y)}{j_{B\cap\Omega'}}\r) \, dy.
\end{equation*}
Using \eqref{eq:key_ineq} and changing variables on the right-hand side, we obtain 
\[
\begin{split}
\int_{\frac{1}{2}B\cap\Omega'}j(y)\log^p\lf(e+\frac{j(y)}{j_{B\cap\Omega'}}\r) \, dy
& \leq C\int_{B\cap\widetilde{\Omega}'}\lf[g(y)\r]^{\frac{np}{n-1}}\lf[j(y)\r]^{1-p}\, dy\\
& \leq C \int_{h^{-1}(B\cap\Omega')}K_h^{\frac{p}{n-1}}(x)\, dx,  
\end{split}
\]
also see~\cite[Lemma 2.2]{KOR:JGA}.

Now we return to the proof of the estimate \eqref{eq:key_ineq} for $B\in\mathcal{B}_1$, 
following an approach similar in spirit to
\cite[Lemma 1.4]{KOR:JGA}. First, we define a distance function by setting 
\begin{equation*}
\mu(y):=\left\{
\begin{array}{ll}
\dist(y, \mathbb R^n\setminus B)  \quad & \text{if}\ y\in \Omega', \\
 0 \quad & \text{otherwise.}
\end{array}
\right.
\end{equation*}
Let $B_r\subset\R^n$ be a ball with radius $r$. We show that
\begin{equation}\label{eq:distanceF}
\bint_{B_r}\mu^n(y)j(y)dy\leq C \left(\bint_{B_{2r}}\mu^{n-1}(y)g(y)dy\right)^{\frac{n}{n-1}}+\widetilde{C}\int_{B\cap\Omega'}j(y)dy,
\end{equation}
where $\widetilde{C}$ depends only on $n$.
Since $\mu(y)=0$ for every $y\in\mathbb R^n\setminus(B\cap\Omega')$, we may assume that the intersection between $B_r$ and $\Omega'$ is not empty. We distinguish between three cases, according to whether $B_{3r}$ is contained in $B\cap\Omega'$ or not.

\textbf{Case $1$:} We assume that $B_{3r}\subset B\cap\Omega'$. By an elementary geometric consideration
we find that
\[\max_{y\in B_r}\mu(y)\leq 4\min_{y\in B_{2r}}\mu(y).\]
Applying (\ref{eq:NEW1}) yields
\begin{eqnarray}
\bint_{B_r}\mu^n(y)j(y)dy&\leq&\max_{y\in B_r}[\mu(y)]^n\bint_{B_r}j(y)dy\nonumber\\
&\leq&C(n)\min_{y\in B_{2r}}[\mu(y)]^n\left(\bint_{B_{2r}}g(y)dy\right)^{\frac{n}{n-1}}\nonumber\\
&\leq&C(n)\left(\bint_{B_{2r}}[\mu(y)]^{n-1}g(y)dy\right)^{\frac{n}{n-1}}.\nonumber
\end{eqnarray}

\textbf{Case $2$:} We assume that $B_{3r}$ is not contained in $B\cap \Omega'$ but is contained in $B$. Recall that $B_r$ intersects $B\cap\Omega'$. An elementary geometric argument shows that
$$
\max_{y\in B_r\cap \Omega'}\mu(y)\leq 4\min_{y\in B_{2r}\cap\Omega'}\mu(y).
$$
Since $\Omega'$ is Ahlfors, $\mu(y)=0$ for $y\notin B\cap\Omega'$ and by \eqref{eq:IIE}, we have
\begin{eqnarray}
\bint_{B_r}\mu^n(y)j(y)dy&\leq&\max_{y\in B_r\cap \Omega'}[\mu(y)]^n\frac{|B_r\cap\Omega'|}{|B_r|}\bint_{B_r\cap\Omega'}j(y)dy\nonumber\\
&\leq&C\min_{y\in B_{2r}\cap \Omega'}[\mu(y)]^n\left(\bint_{B_{2r}\cap\Omega'}g(y)dy\right)^{\frac{n}{n-1}}\nonumber\\
&\leq&C\left(\bint_{B_{2r}}[\mu(y)]^{n-1}g(y)dy\right)^{\frac{n}{n-1}}.\nonumber
\end{eqnarray}
Here the constant $C$ depends on all the parameters named in Lemma \ref{le:IIE}.

\textbf{Case $3$:} We assume that $B_{3r}$ is not contained in $B$. Again, recall that $B_r$ intersects $B\cap\Omega'$. An elementary geometric argument shows that
$$
\max_{y\in B_r}\mu(y)\leq 4 r.
$$ 
Since $\mu(y)=0$ for $y\notin B\cap\Omega'$, we have
\begin{eqnarray}
\bint_{B_r}[\mu(y)]^nj(y)dy&\leq&\max_{y\in B_r}[\mu(y)]^n\frac{1}{|B_r|}\int_{B_r\cap B\cap\Omega'}j(y)dy\nonumber\\
&\leq&C(n)\frac{r^n}{|B_r|}\int_{B\cap\Omega'}j(y)dy\nonumber\\
&\leq& C(n)\int_{B\cap\Omega'}j(y)dy.\nonumber
\end{eqnarray}
Combining these three cases proves the inequality (\ref{eq:distanceF}).

We split our argument into several steps now. 

\textbf{Step 1: Using the maximal operator.} 
For $y\in\mathbb{R}^n$, the estimate \eqref{eq:distanceF} implies (after taking the supremum over all $r$
and raising the obtained inequality to the power $1/n$) the following point-wise
estimate 
\begin{equation}\label{eq:maximal}
\left[\mathbf{M}[\mu^n j](y)\right]^{\frac{1}{n}}\leq C\left[\mathbf{M}[\mu^{n-1}g](y)\right]^{\frac{1}{n-1}}+\widetilde{C}\left(\int_{B\cap\Omega'}j(y)dy\right)^{\frac{1}{n}},
\end{equation}
where $\mathbf M$ stands for the classical Hardy-Littlewood maximal operator. Recall that we assumed \eqref{eq:normalization} to simplify further calculations, i.e., the second term on the right-hand side is without loss of generality equal to $1/2$.
Taking this normalization into account, and applying \eqref{eq:maximal}, for $t>1$ we have 
\begin{equation}\label{eq:maximal1}
\abs{\lf\{ y\in \mathbb R^n:\mathbf{M}[\mu^nj](y)>t^n\r\}}
\leq \abs{\lf\{y\in \mathbb R^n: C\mathbf M[\mu^{n-1}g](y)>t^{n-1}\r\}}.
\end{equation}

As a classical application of Vitali's covering lemma (see e.g.\,\cite{stein1970}), for every $h\in L^1(\mathbb R^n)$ and $\lambda>0$ we have 
\begin{equation}\label{eq:weak1-1}
\left|\left\{x\in\mathbb R^n: \mathbf{M}[h](x)>\lambda\right\}\right|\leq\frac{C(n)}{\lambda}\int_{\{x\in\mathbb R^n:2|h(x)|>\lambda\}}|h(y)|dy.
\end{equation} 
As a direct consequence of the Calder\'{o}n-Zygmund decomposition (see e.g.\,\cite{stein1970}), the reverse inequality also holds. Precisely, for every $h\in L^1(\mathbb R^n)$ and $\lambda>0$ we have 
\begin{equation}\label{Rweak1-1}
C(n)\left|\left\{x\in\mathbb R^n:\mathbf{M}[h](x)>\lambda\right\}\right|\geq\frac{1}{\lambda}\int_{\left\{x\in\mathbb R^n:|h(x)|>\lambda\right\}}|h(y)|dx.
\end{equation}
By applying inequalities (\ref{eq:weak1-1}) and (\ref{Rweak1-1}) to (\ref{eq:maximal1}), we obtain that there exists a constant $\gamma>1$,  depending only on $n$, such that for every $t>1$ 
\begin{equation}\label{eq:11bounded}
    \int_{\lf\{y\in \mathbb{R}^n:\mu^n(y)j(y)>t^n\r\}}\mu^n(y)j(y)dy
    \leq C t\int_{\lf\{y\in \mathbb{R}^n:\gamma\mu^{n-1}(y)g(y)>t^{n-1}\r\}}\mu^{n-1}(y)g(y)dy.
\end{equation}

\textbf{Step 2: Using an auxiliary logarithmic function.} 
Let $1\leq p<\infty$. Define an auxiliary function by setting 
\[\Phi(t):=\frac{1}{p}\log^pt+\log^{p-1}t,\ \ t\geq 1.\]
Then for $t>1$ we have 
\begin{equation}\label{eq:equa2}
\Phi'(t)>0\ \ {\rm and}\ \ \frac{d}{dt}\lf(t\log^{p-1}t\r)=t\Phi'(t).
\end{equation}
We multiply both sides of the inequality (\ref{eq:11bounded}) by $\Phi'(t)$, fix $T>2$ and then integrate the resulting inequality over $[1,T]$.
Applying (\ref{eq:equa2}) to the left-hand side and using the fact that $\mu(y)=0$ if $y\notin B\cap\Omega'$, Fubini's theorem yields 
\begin{equation}\label{eq:equa1}
\begin{split}
&\int_{\{y\in B\cap\Omega':1<[\mu(y)]^nj(y)<T^n\}}\mu^n(y)j(y)\cdot\left(\frac{1}{p}\log^p(\mu(y) [j(y)]^{1/n})\right)dy\\ 
&\quad\le \int_1^T\Phi'(t)\int_{\lf\{y\in \mathbb{R}^n:\mu^n(y)j(y)>t^n\r\}}\mu^n(y)j(y)dy \, dt\\
&\quad\leq C\int_1^T \Phi'(t)t\int_{\lf\{y\in \mathbb{R}^n:\gamma[\mu(y)]^{n-1}g(y)>t^{n-1}\r\}}\mu^{n-1}(y)g(y)dy\, dt\\
&\quad\le C\int_{\lf\{y\in B\cap\Omega': 1<\gamma[\mu(y)]^{n-1}g(y)<T^{n-1}\r\}}[\mu(y)]^n [g(y)]^{\frac{n}{n-1}}\log^{p-1}\lf(e+\gamma^{\frac{1}{n-1}}\mu(y) [g(y)]^{\frac{1}{n-1}}\r).
\end{split}
\end{equation}
Note that $j(y)=0$ necessarily implies that $g(y)=0$. Moreover, for every $y$ we have $g(y)\geq j(y)^{\frac{n-1}{n}}$. Hence, 
\[\lf\{y\in B\cap\Omega': 1<\gamma[\mu(y)]^{n-1}g(y)<T^{n-1}\r\}\]
is a subset of 
\[E:=\lf\{y\in B\cap\Omega': T^n> [\mu(y)]^nj(y)\ {\rm and}\ j(y)>0\r\}.\]
Then we have 
\begin{equation}
\int_{\lf\{y\in B\cap\Omega': 1<\gamma[\mu(y)]^{n-1}g(y)<T^{n-1}\r\}}[\mu(y)]^n[g(y)]^{\frac{n}{n-1}}\log^{p-1}\lf(e+\gamma^{\frac{1}{n-1}}\mu(y) [g(y)]^{\frac{1}{n-1}}\r)dy
\leq C(I_1+I_2)\nonumber
\end{equation}
with 
\[I_1:=\int_E[\mu(y)]^n[g(y)]^{\frac{n}{n-1}}\log^{p-1}\lf(e+\mu(y)[j(y)]^{\frac{1}{n}}\r)dy\]
and
\[I_2:=\int_E[\mu(y)]^n[g(y)]^{\frac{n}{n-1}}\log^{p-1}\lf(e+[\gamma g(y)]^{\frac{1}{n-1}}[j(y)]^{-\frac{1}{n}}\r)dy.\]
Since $g^{\frac{n}{n-1}}(y)\geq j(y)$,  there exists a constant $C>1$ independent of $p$ with
\begin{equation}
\int_{\lf\{y\in B\cap\Omega': [\mu(y)]^nj(y)< 2 \r\}}[\mu(y)]^nj(y)\log^p\lf(e+\mu(y)[j(y)]^{\frac{1}{n}}\r)dy\leq C I_1.
\end{equation}
Applying (\ref{eq:equa1}) gives
\begin{equation}\label{eq:equaI12}
    \int_{\{y\in B\cap\Omega':[\mu(y)]^nj(y)<T^n\}}[\mu(y)]^nj(y)\log^p\lf(e+\mu(y)[j(y)]^{\frac{1}{n}}\r)dy\\
    \leq C\lf(I_1+I_2\r).
\end{equation}

\textbf{Step 3: Estimating $I_1$ and $I_2$.} 
Here we use $g(y)\geq[j(y)]^{\frac{n-1}{n}}$ (making the argument of the logarithm in $I_2$ at least $e+1$) and 
$[\mu(y)]^n\leq C \abs{B\cap\Omega'}\chi_{B\cap\Omega'}(y)$. 
For $I_1$, Young's inequality yields
\begin{equation}\label{eq:I1}
\begin{split}
I_1& = \int_E\left[[\mu(y)]^{\frac{n}{p}}[g(y)]^{\frac{n}{n-1}}j(y)^{-\frac{p-1}{p}} \right]\cdot\left[[\mu(y)]^{\frac{n(p-1)}{p}}j(y)^{\frac{p-1}{p}}\log^{p-1}\lf(e+\mu(y)[j(y)]^{\frac{1}{n}}\r)\right]dy\\
&\leq\frac{  C\abs{B\cap\Omega'}}{\varepsilon^p}\int_{B\cap\widetilde{\Omega'}}[g(y)]^{\frac{np}{n-1}}[j(y)]^{1-p}dy +C\varepsilon^{\frac{p}{p-1}}\int_E[\mu(y)]^nj(y)\log^p\lf(e+\mu(y)[j(y)]^{\frac{1}{n}}\r)dy,
\end{split}
\end{equation} 
where $\varepsilon$ denotes a small positive constant that will be fixed later.
For $I_2$, we have 
\begin{equation}\label{eq:I2}
    I_2\leq C |B\cap\Omega'|\int_{B\cap\widetilde{\Omega'}}[g(y)]^{\frac{np}{n-1}}[j(y)]^{1-p}dy.
\end{equation}

Now, we choose $\varepsilon$ to be sufficiently small such that the second term on the right-hand side of (\ref{eq:I1}) can be absorbed into the term on the left-hand side of (\ref{eq:equaI12}) (this can be done as  $\varepsilon$ then depends only on the same parameters as $C$). 
Letting $T\rightarrow \infty$ and applying the monotone convergence theorem, we obtain 
\begin{equation}\label{eq:equa4}
\begin{split}
    &\int_{B\cap\Omega'}[\mu(y)]^nj(y)\log^p\lf(e+\mu(y)[j(y)]^{\frac{1}{n}}\r)dy\\
    &\quad\leq C  \abs{B\cap\Omega'}\int_{B\cap\widetilde{\Omega'}}[g(y)]^{\frac{np}{n-1}}[j(y)]^{1-p}dy.
\end{split}
\end{equation}
Hence, we have 
\begin{equation}\label{eq:equa5}
\begin{split}
&\int_{\frac{1}{2}B\cap\Omega'}[\mu(y)]^nj(y)\log^p\lf(e+  \abs{B\cap\Omega'}[j(y)]^{\frac{1}{n}}\r)dy\\
&\quad\leq C  \abs{B\cap\Omega'} \int_{B\cap\widetilde{\Omega'}}[g(y)]^{\frac{np}{n-1}}[j(y)]^{1-p}dy.
\end{split}
\end{equation}
This is because for every $y\in\frac{1}{2}B\cap\Omega'$, we have 
\[[\mu(y)]^n\geq\lf(\frac{1}{2}r_B\r)^n\geq C \abs{B\cap\Omega'},\]
where $r_B>0$ denotes the radius of the ball $B$.

Hence,  
\begin{align*}
\int_{\frac{1}{2}B\cap\Omega'}j(y)\log^p\lf(e+  \abs{B\cap\Omega'}[j(y)]^{\frac{1}{n}}\r)dy
    \leq C\int_{B\cap\widetilde{\Omega'}}[g(y)]^{\frac{np}{n-1}}[j(y)]^{1-p}dy.\nonumber
\end{align*}
Finally, taking the normalization (\ref{eq:normalization}) into account, we obtain the desired result
\begin{equation*}
\begin{aligned}
\int_{\frac{1}{2}B\cap\Omega'}j(y)\log^p\lf(e+\frac{j(y)}{j_{B\cap\Omega'}}\r)dy &=\int_{\frac{1}{2}B\cap\Omega'}j(y)\log^p\lf(e+C|B\cap\Omega'|j(y)\r)dy\\
&\le C\int_{\frac{1}{2}B\cap\Omega'}j(y)\log^p\lf(e+|B\cap\Omega'|j^{\frac{1}{n}}(y)\r)dy\\
&\leq C\int_{B\cap\widetilde\Omega'}\lf[g(y)\r]^{\frac{np}{n-1}}\lf[j(y)\r]^{1-p}dy.\\
\end{aligned}
\end{equation*}

\end{proof}

Given the finite cardinality of 
$\mathcal B$, there exists a sufficiently large constant $C$ ensuring that inequality~(\ref{eq:ball}) holds with a constant independent of the ball $B$. Now, we are ready to prove Theorem~\ref{thm:intejaco_main}. 
\begin{proof}[Proof of Theorem \ref{thm:intejaco_main}]
Since
\[\overline{\Omega'}\subset\bigcup_{B\in\mathcal B}\frac{1}{2}B\]
 and 
 $h:\overline{\Omega}\to\overline{\Omega'}$ is a homeomorphism, we have 
\[\Omega\subset\bigcup_{B\in\mathcal B}h^{-1}\lf(\frac{1}{2}B\r).\]
Then, by summing inequalities (\ref{eq:ball}) over balls in $\mathcal B$, we have
\begin{eqnarray}
\int_\Omega\log^p\lf(e+\frac{1}{J_h(x)}\r)dx&\leq&\sum_{B\in\mathcal B}\int_{h^{-1}\lf(\frac{1}{2}B\cap\Omega'\r)}\log^p\lf(e+\frac{1}{J_h(x)}\r)dx\nonumber\\
&\leq&C\sum_{B\in\mathcal B}\int_{h^{-1}\lf(B\cap\Omega'\r)}K^{\frac{p}{n-1}}_h(x)dx\nonumber\\
&\leq&C\int_\Omega K^{\frac{p}{n-1}}_h(x)dx.\nonumber
\end{eqnarray}
\end{proof}

\section{Proof of Theorem~\ref{thm:outdistortion}}
Let $\phi$ be a positive convex function on $(0, \infty)$ with
\begin{equation}\label{eq:condition1}
    \lim_{t\to0^+}\phi(t)=\infty\ {\rm and}\ \phi(t)=\infty\ {\rm for}\ t\leq 0.
\end{equation}
Also assume there is $A>0$ such that
\begin{equation}\label{eq:condition2}
    A^{-1}\phi(t)\leq\phi(2t)\leq A\phi(t),\ {\rm for}\ t\in(0, \infty).
\end{equation}
Then for every homeomorphism $h:\Omega\to\Omega'$ between domains $\Omega,\Omega'\subset\rn$ in the class $W^{1, n-1}(\Omega, \Omega')$, we have the corresponding energy 
\begin{equation}\label{eq:energy}
   \mathcal E^{n-1}(h):=\int_{\Omega}\lf(|Dh(x)|^{n-1}+\phi(J_h(x))\r)dx, 
\end{equation}
see \eqref{eq:neo}.
The following result is from \cite{INV}.
\begin{thm}
\label{prop:}
Let $\Omega, \Omega'\subset\rn$ be Lipschitz domains with $n\geq 3$. Let $\phi$ be a positive convex function on $(0, \infty)$ which satisfies (\ref{eq:condition1}) and (\ref{eq:condition2}). Let $\{h_k\}_{k=1}^{\infty}\subset W^{1, n-1}(\Omega, \Omega')$ be a sequence of homeomorphisms with $J_{h_k}>0$ almost everywhere and 
\[\sup_k\mathcal E^{n-1}(h_k)<\infty\ {\rm and}\ \sup_{k}\int_{\Omega}K_{h_k}^{\frac{1}{n-1}}(x)dx<\infty.\]
Let $\{h_k\}$ converge weakly in $W^{1, n-1}(\Omega, \Omega')$ to a limit mapping $h$. 
Assume further that  
 $h_k$ are homeomorphisms of $\overline\Omega$ onto $\overline{\Omega'}$ such that $h_k=h_1$ on $\partial\Omega$ for all $k\in\mathbb N$.
Then $h$ satisfies the $(INV)$ condition.
\end{thm}

Theorem \ref{thm:outdistortion} improves this result in the following sense: to gain the $(INV)$ condition of the weak limit $h$, we do not need to assume the uniform boundedness of the energy term. Theorem \ref{thm:intejaco_main} shows that for the Jacobian term it is implied by the integrability of distortion functions $K_{h_k}$, and the derivative term is bounded from the weak convergence. 

\begin{proof}[Proof of Theorem \ref{thm:outdistortion}]
Define 
\[\phi(t):=\begin{cases}
    \log\lf(e+\frac{1}{t}\r), &\ {\rm if}\ t\in (0, \infty),\\
     \infty, &\ t\in(-\infty, 0].
     \end{cases}\]
One can verify that $\phi$ is a positive convex function on $(0, \infty)$ which satisfies (\ref{eq:condition1}) and (\ref{eq:condition2}). Let $\{h_k\}_{k=1}^\infty$ be a sequence of homeomorphisms which satisfies the assumption of Theorem \ref{thm:outdistortion}. By Theorem \ref{thm:intejaco_main}, we have 
\[\sup_{k}\int_\Omega\phi(J_{h_k}(x))dx<\infty.\]
Since $h_k$ converge weakly to $h$ in $W^{1, n-1}(\Omega, \Omega')$, we have 
\[\sup_{k}\int_{\Omega}|Dh_k(x)|^{n-1}dx<\infty.\]
Hence, the sequence $\{h_k\}$ satisfies the assumptions in Theorem \ref{prop:}, so $h$ satisfies the $(INV)$ condition.
\end{proof}

The discussion in next remark shows the sharpness of $L^{\frac{1}{n-1}}$-integrability of distortion functions.

\begin{remark}\label{rk:integrabilty_K}
In \cite{INV}, for every $0<\alpha<2$, there is a sequence of homeomorphisms $\{h_k\}_{k=1}^\infty$ from $\overline{B^3(0, 10)}$ onto $\overline{B^3(0, 10)}$ such that 
$h_k\in W^{1, 2}(B^3(0, 10), B^3(0, 10))$ with $h_k=id$ on the boundary $\partial B^3(0, 10)$ for every $k\in\mathbb N$ and 
\[\sup_k\int_{B^3(0, 10)}|Dh_k(x)|^2+\frac{1}{J_{h_k}^\alpha(x)}dx<\infty.\]
However, the weak limit of $\{h_k\}$ does not satisfy the $(INV)$ condition. Take $0<\beta<1/2$ and choose $\alpha=\frac{2\beta}{2-3\beta}$, then the H\"older inequality implies that 
\begin{eqnarray}
\int_{B^3(0, 10)}K^\beta_{h_k}(x)dx&=&\int_{B^3(0, 10)}\frac{|Dh_k(x)|^{3\beta}}{J^\beta_{h_k}(x)}dx\nonumber\\
                                                 &\leq&\lf(\int_{B^3(0, 10)}|Dh_k(x)|^2dx\r)^{\frac{3\beta}{2}}\lf(\int_{B^3(0, 10)}\lf(\frac{1}{J_{h_k}(x)}\r)^{\frac{2\beta}{2-3\beta}}dx\r)^{\frac{2-3\beta}{2}}\nonumber\\
                                                   &<&\infty.\nonumber
\end{eqnarray}
As $\alpha<2$, the weak limit of $\{h_k\}$ does not satisfy the $(INV)$ condition, it shows that (at least for $n=3$) one cannot lower the exponent in (\ref{eq:supremum}).
\end{remark}

\end{document}